\documentclass[11pt]{article}

\usepackage{amsmath}
\usepackage{amsopn}
\usepackage{latexsym}
\usepackage{amsfonts}
\usepackage{amssymb}
\usepackage{graphicx}
\usepackage{subcaption}
\usepackage{amsthm}

\usepackage{tikz}
\usetikzlibrary{shapes.geometric, intersections}
\usepackage{mathdots}
\usepackage{mathtools}
\usepackage{bm}
\usepackage{xcolor}
\usepackage{enumitem}

\DeclareMathOperator*{\argmin}{arg\,min}

\allowdisplaybreaks

\newtheorem{claim}{Claim}
\newtheorem{definition}{Definition}
\newtheorem*{remark}{Remark}

\title{High Order Schemes for \\ Gradient Flow with Respect to a Metric}
\author{Saem Han\footnote{Department of Mathematics, University of Michigan, Ann Arbor, MI 48109, USA; saem@umich.edu (corresponding author)} \and Selim Esedo\=glu\footnote{Department of Mathematics, University of Michigan, Ann Arbor, MI 48109, USA; esedoglu@umich.edu} \and Krishna Garikipati\footnote{Departments of Mechanical Engineering, and Mathematics, Michigan Institute for Computational Discovery \& Engineering, University of Michigan, Ann Arbor, MI 48109, USA; krishna@umich.edu}}
\date{}

\begin{document}
\maketitle

\begin{abstract}
New criteria for energy stability of multi-step, multi-stage, and mixed schemes are introduced in the context of evolution equations that arise as gradient flow with respect to a metric.
These criteria are used to exhibit second and third order consistent, energy stable schemes, which are then demonstrated on several partial differential equations that arise as gradient flow with respect to the 2-Wasserstein metric.
\end{abstract}


\section{Introduction}
\label{sec:intro}

Let a cost function, or energy, $\mathcal{E} \, : \, X \to \mathbb{R}$ be defined on a metric space $(X,d)$.
One way to define gradient flow of $\mathcal{E}$ on $X$ with respect to the metric $d$ is as the $k\to 0^+$ limit of discrete in time approximations obtained by solving the sequence of variational problems
\begin{equation}
\label{eq:minmov}
u_{n+1} \in \argmin_{\xi\in X} \left\{ \mathcal{E}(\xi) + \frac{1}{2k} d^2(\xi,u_n)  \right\}
\end{equation}
starting from an initial condition $u_0\in X$.
This approach to extending theory and numerics for stationary problems involving an energy $\mathcal{E}$ to dynamics is often referred to as {\em minimizing movements} \cite{degiorgi, atw, luckhaus}.
Conditions are needed on $(X,d)$ and $\mathcal{E}$ to ensure well-posedness of the optimizations (\ref{eq:minmov}), which should hold at least for small enough time steps $k>0$.

In the most standard setting where $X$ is an inner product space so that $d(\xi,\eta) = \langle  \xi - \eta \, , \, \xi - \eta \rangle^{\frac{1}{2}}$, equation (\ref{eq:minmov}) is the variational formulation of the standard backward-Euler discretization
\begin{equation}
u_{n+1} = u_n - k \nabla_X \mathcal{E}(u_{n+1})
\end{equation}
of the abstract ODE
\begin{equation}
\label{eq:gradientFlow}
u'(t) = -\nabla_X \mathcal{E}(u)
\end{equation}
the standard definition of which requires certain differentiability conditions on $\mathcal{E}$.

Assuming the $k\to 0^+$ limit of the time discrete solutions given by (\ref{eq:minmov}) exists, and taking that limit as the exact solution, in this paper we study related variational schemes that aim to converge to the same limit, but {\em faster}.
They can be seen as multi-step and/or multi-stage analogues of (\ref{eq:minmov}), and take the form
\begin{equation}
\label{eq:step}
u_{n+1} \in \argmin_{\xi\in X} \left\{ \mathcal{E}(\xi) + \frac{1}{k} \sum_{j=1}^K \alpha_j d^2(\xi,\eta_j)\right\}
\end{equation}
for some $\eta_j\in X$ ultimately depending only on $u_0,u_1,\ldots,u_n$, and constants $\alpha_j\in\mathbb{R}$.

We will assume that the stationary problem (\ref{eq:step}) has a unique solution for $k>0$ small enough, and that there is already an effective optimization algorithm for finding it.
Here, the focus of our investigation is {\em stability} and {\em consistency}:
Are there conditions on the coefficients $\alpha_j$ that ensure dissipation or equiboundedness as $k\to 0^+$ of the energy, as well as high order consistency whenever the solution is sufficiently smooth, merely by virtue of $d$ being a metric?

Our contribution is thus simple, sufficient criteria  to ensure energy dissipation and energy boundedness of multi-step, multi-stage, and mixed schemes that utilize (\ref{eq:step}) as their sole building block.
The criteria are based on embedding certain finite graphs associated with the given scheme into Euclidean spaces.
As a concrete application, we use these criteria to search for and exhibit a new, energy decreasing multistage (Runge-Kutta) scheme that is second order accurate, and a new mixed (multi-step / multi-stage) scheme that is energy bounded and third order accurate, for any metric $d$ in (\ref{eq:step}).
Taking the 2-Wasserstein distance on the space of probability measures as our example of a metric space, we demonstrate the new schemes on the heat, porous medium, and nonlinear Fokker-Planck equations.
Our work can be seen as a higher order follow up to the investigation in \cite{matthes_plazotta} that studies the second order BDF2 scheme in the general context of metric spaces (see also the very recent contribution \cite{gallouet} that studies BDF2, and hence second order schemes exclusively, unlike this paper).
We refrain, however, from a discussion of general conditions that ensure uniqueness of minimizers of (\ref{eq:step}) and leave it to be checked on a problem by problem basis, while acknowleding the question can be quite involved, as discussed in \cite{matthes_plazotta}, especially without imposing additional conditions on the time step size $k$.

\section{Previous Work}
\label{sec:previous}
There are numerous studies in existing literature devoted to energy stability of implicit as well as semi-implicit schemes for evolution problems that arise as gradient descent for an energy, formulated in the standard context of gradient flow with respect to an inner product.
Closest to the approach of the present paper are those that can be seen as multi-stage or multi-step generalizations of the variational formulation of backward Euler, or minimizing movements, scheme (\ref{eq:minmov}).
We note the early contribution \cite{hairer_lubich} as a starting point.

Recently, there has been interest in high order accurate, energy stable schemes that can be applied universally for gradient flows in any metric.
Particular attention has been paid to Wassertein gradient flows and their connection to important partial differential equations, e.g. Fokker-Planck and porous medium equations, among others, following the pioneering work of \cite{jko}.
To date, these schemes are at most second order accurate in time, and typically only guarantee boundedness (rather than dissipation) of the energy.
We note three papers in particular, namely \cite{legendre_turinici, matthes_plazotta, carillo_craig_wang_wei}, that are closely related to the present contribution.

In \cite{legendre_turinici}, a variational formulation of the standard midpoint formula is discussed, and is applied to Wasserstein gradient flows, resulting in a second order consistent and unconditionally stable (energy bounded) numerical scheme.
The approach of \cite{matthes_plazotta} is based on the standard backward differentiation formula BDF2, which is shown to have a variational formulation (a multi-step version of minimizing movements) when applied to gradient flows.
It yields a second order consistent scheme that once again enjoys, in terms of stability, boundedness of its energy.
Finally, \cite{carillo_craig_wang_wei} considers adapting the standard Crank-Nicholson scheme to the context of Wasserstein gradient descent.
However, it is not clear that the resulting scheme achieves better than first order accuracy in time; the focus of the authors is rather on efficient solution of the variational problem involved at every time step.

The present work is closest to \cite{matthes_plazotta}, where authors point out their particular approach based on backward differentiation formulas will not yield third and higher order schemes that enjoy the same desirable stability properties.
It is then natural to ask whether multi-stage or mixed (multi-stage / multi-step) analogues that rely on similar optimization problem at every time step as the scheme of \cite{matthes_plazotta}, can be found that achieve third order consistency as well as energy stability.
After developing a new approach to energy stability in Section \ref{sec:stability}, we exhibit such schemes in Section \ref{sec:consistency}.

\section{Stability}
\label{sec:stability}
In this section, we introduce a general mixed scheme for gradient flows which includes both multi-step schemes such as backward differentiation formulas and multi-stage schemes such as Runge-Kutta methods as special cases, and study conditions for this scheme to be stable. Let $(X,d)$ be a metric space and $\mathcal{E}:X\to (-\infty,\infty]$ be an energy functional which is bounded from below. The mixed scheme for gradient flow (\ref{eq:gradientFlow}) using $M$ previous steps and $N$ intermediate stages, $M,N\in\mathbb{N}$, is formulated as follows. Starting from the initial time $t=t_0$, for a fixed time step size $k$, we denote by $u_j$ the discrete solution of (\ref{eq:gradientFlow}) at the $j$-th time step: $u_j = u(\cdot, t_0+jk)$. Given $M$ previous steps $u_{n-M+1},\cdots,u_{n}$, we find the next step $u_{n+1}$ by constructing $N$ intermediate stages using up to $M$ previous steps. The $N$ intermediate stages are obtained by solving optimization problems that minimize the original energy functional $\mathcal{E}$ plus additional terms which we call movement limiters.
For notational convenience, we will use non-positive indices to denote previous steps. 
\begin{align}
    &\text{For } i=-M+1,\cdots,0: \nonumber\\
    &\quad \text{ Set } v_i := u_{n+i}. \nonumber\\
    &\text{For } i=1,\cdots,N: \nonumber\\
    &\quad v_i = \argmin_{\xi\in X} \left\{ \mathcal{E}(\xi) + \frac{1}{2k}\sum_{j=-M+1}^{i-1} \gamma_{i,j} d^2(\xi,v_j)\right\}.\nonumber \\
    &\text{Set } u_{n+1}:=v_N 
    \label{eq:scheme}
\end{align}
for $\gamma_{i,j}\in\mathbb{R}$, $i=1,\cdots,N$, $j=-M+1,\cdots,i-1$. 
Of course, in general, additional conditions need to be imposed on $X$, $d$, $\mathcal{E}$, or $k$ to ensure variational problems in (\ref{eq:scheme}) are uniquely solvable.
For the present, we assume this has been verified, and investigate the conditions on $\gamma_{i,j}$'s which ensure that scheme (\ref{eq:scheme}) preserves the principal qualitative feature of gradient flows, namely energy dissipation, which is often referred to as energy stability in the context of numerical schemes:
\begin{equation}
    \mathcal{E}(u_{n+1}) \leq \mathcal{E}(u_n).
    \label{eq:energyStability}
\end{equation}
Short of (\ref{eq:energyStability}), we will subsequently also consider the less stringent stability property of an energy bound that holds uniformly for all small enough time steps $k$.

From the fact that each intermediate stage $v_i$ minimizes the energy $\mathcal{E}$ plus movement limiter terms, we have the following inequality for $i=1,\cdots,N$:
\begin{equation}
    \mathcal{E}(v_i)+\frac{1}{2k}\sum_{j=-M+1}^{i-1}\gamma_{i,j}d^2(v_i,v_j)\leq \mathcal{E}(v_{i-1})+\frac{1}{2k}\sum_{j=-M+1}^{i-2}\gamma_{i,j}d^2(v_{i-1},v_j).
\end{equation}
By induction,
\begin{align}
    \mathcal{E}(v_N)+\frac{1}{2k}\sum_{i=1}^N \sum_{j=-M+1}^{i-1}\gamma_{i,j}d^2(v_i,v_j) \leq \mathcal{E}(v_0)+\frac{1}{2k}\sum_{i=1}^N \sum_{j=-M+1}^{i-2} \gamma_{i,j} d^2(v_{i-1},v_j),
\end{align}
and by rearranging terms, we obtain
\begin{align}
    \mathcal{E}(v_N)\leq \mathcal{E}(v_0)+\frac{1}{2k}\left\{\sum_{i=0}^{N-1}\sum_{j=-M+1}^{i-1}\gamma_{i+1,j}d^2(v_i,v_j)-\sum_{i=1}^N\sum_{j=-M+1}^{i-1}\gamma_{i,j}d^2(v_i,v_j)\right\}.
\end{align}
Define 
\begin{equation}
    w_{i,j}=\begin{dcases}
        \gamma_{i+1,j}-\gamma_{i,j}, \quad &\text{for } i>j\\
        0, &\text{for } i\leq j
    \end{dcases}
    \label{eq:weights}
\end{equation}
for $i=0,\cdots,N$ and $j=-M+1,\cdots,N$ by assuming $\gamma_{0,j}=0$ and $\gamma_{N+1,j}=0$ for all $j$. Since $v_0=u_n$ and $v_N=u_{n+1}$, it follows that 
\begin{equation}
    \mathcal{E}(u_{n+1})\leq \mathcal{E}(u_n)+\frac{1}{2k} \sum_{i=0}^{N} \sum_{j=-M+1}^{N} w_{i,j} d^2(v_{i},v_{j}).
    \label{eq:stabilityIneq}
\end{equation}
Therefore, the scheme (\ref{eq:scheme}) is energy stable provided that
\begin{equation}
\sum_{i=0}^{N} \sum_{j=-M+1}^{N} w_{i,j} d^2(v_{i},v_{j})\leq 0.
\end{equation}

We aim to find as weak conditions on $\gamma_{i,j}$ as possible that ensure stability. In case the metric $d$ arises from an inner product, i.e. $d^2(p,q) = \langle p - q \, , \, p - q \rangle$, for any $\{ v_j \}_{j=-M+1}^N \subset X$ there exists $\{ p_j \}_{j=-M+1}^N \subset \mathbb{R}^m$ for some $m$ so that $|p_i - p_j| = d(v_i,v_j)$ for all $i,j$.
Then, the additional terms
\begin{equation}
\label{eq:additional}
\sum_{i=0}^{N} \sum_{j=-M+1}^{N} w_{i,j} d^2(v_{i},v_{j})
\end{equation}
in formula (\ref{eq:stabilityIneq}) are given by the quadratic form
\begin{equation}
\sum_k \sum_{i,j} w_{i,j} \left( (p_i)_k - (p_j)_k \right)^2.
\end{equation}
In that case, the sign of (\ref{eq:additional}) can be investigated via the eigenvalues of the matrix
\begin{equation}
\sum_{i,j} w_{i,j} \left( e_{i+M} - e_{j+M} \right) \otimes \left( e_{i+M} - e_{j+M} \right)
\end{equation}
representing the quadratic form; here $e_i$ denotes the $i$-th standard basis element of $\mathbb{R}^{M+N}$.
For instance, energy dissipation would be ensured if the matrix is negative semi-definite.
This is the gist of the approach in e.g. \cite{hairer_lubich}.

For a general metric $d$, all we know are {\bf 1.} the terms $d(v_i,v_j)$ are non-negative, and {\bf 2.} they satisfy the triangle inequality.
Given the complete, undirected graph for which the pairwise distances $d_{i,j} = d(v_i,v_j)$ form the edge weights between the $i$-th and $j$-th vertices, we cannot in general find $p_j$ in $\mathbb{R}^m$ for {\em any} $m$ such that $|p_i-p_j| = d_{i,j}$: Not all finite metric spaces are isometrically embeddable into a Euclidean space.
However, we may still be able to bound (\ref{eq:additional}) from above by a quadratic form: We can partition (\ref{eq:additional}) into terms each of which corresponds to a subgraph, obtained by removing some of the edges, that is always embeddable.

\begin{definition}
For us, an undirected graph $G$ on a set of $N$ vertices is a pair $G=(V,E)$ where $V\subset \mathbb{Z}$ with $\# V = N$, $E\subset V\times V$ such that $(i,j)\in E \rightarrow (j,i)\in E$ and $(i,i)\not\in E$ for any $i=1,2,\cdots,N$. 
\end{definition}

\begin{definition}
We will say that an undirected graph $G=(V,E)$ on a set of $N$ vertices is embaddable in $\mathbb{R}^m$ for some $m$ if given any set of pairwise distances $\{d_{i,j}\}_{(i,j)\in E}$ satisfying
\begin{enumerate}
    \item $d_{i,j} \geq 0$,
    \item $d_{i,j}+d_{j,k}\geq d_{i,k}$ whenever $(i,j), (j,k),(i,k)
    \in E$ and $i,j,k$ are distinct,
\end{enumerate}
there exists $p_1,\cdots,p_N\in \mathbb{R}^m$ such that
\begin{equation}
    |p_i-p_j|=d_{i,j} \quad \text{for each } (i,j)\in E.
\end{equation}
\end{definition}

\textbf{Example.} The complete graph on four vertices is not embaddable into any $\mathbb{R}^m$ but the graph on $4$ vertices with $V=\{1,2,3,4\}$ and 
\begin{equation}
    E=\{(1,2),(2,1),(1,3),(3,1),(2,3),(3,2),(3,4),(4,3),(2,4),(4,2)\}
\end{equation}
is embeddable in $\mathbb{R}^2$.

\begin{remark}
    If $G$ is embaddable, then so is any of its subgraphs.
\end{remark}

The following claim is simply an outline of how we intend to approach the stability question of the general mixed (multi-step /  multi-stage) scheme of the form (\ref{eq:scheme}).

\begin{claim}
\label{claim:stability}
Let $(X,d)$ be a metric space. Given a set of coefficients $\{\gamma_{i,j}\in \mathbb{R}:i=1,\cdots,N$, $j=-M+1,\cdots,i-1\}$, define $w_{i,j}\in \mathbb{R}$ for $i=0,\cdots,N$ and $j=-M+1,\cdots,N$ as in (\ref{eq:weights}). 

Let $G=(V,E)$ be an undirected graph on a set of $(M+N)$ vertices with $V=\{-M+1,\cdots,N\}$ and
\begin{equation}
    E=\{(i,j)\subset V\times V: \text{either } w_{i,j}\neq 0 \text{ or } w_{j,i}\neq 0\}.
\end{equation}
Consider a collection of subgraphs $\{S^\alpha\}$ of $G$ such that 
\begin{equation}
    S^\alpha = (V,E^\alpha) \text{ and } E = \cup_\alpha E^\alpha,
\end{equation}
i.e., every edge of $G$ belongs to at least one $S^\alpha$. 

For each $\alpha$, choose $\theta^\alpha_{i,j}\in \mathbb{R}$ for $(i,j)\in E^\alpha$ satisfying  
\begin{equation}
    \sum_{\alpha} \theta_{i,j}^\alpha=w_{i,j},
    \label{eq:weightPartition}
\end{equation}
and define a symmetric matrix $Q^\alpha$ as follows: 
\begin{equation}
    Q^\alpha = \sum_{(i,j)\in E^\alpha} \theta_{i,j}^\alpha (e_{i+M}-e_{j+M})\otimes (e_{i+M}-e_{j+M}).
\end{equation}
Suppose that the collection of subgraphs $\{S^\alpha\}$ and the corresponding set of matrices $\{Q^\alpha\}$ satisfy the following:
\begin{enumerate}[label=Assumption \arabic*., leftmargin=*]
    \item All $S^\alpha$'s are embeddable into $\mathbb{R}^m$ for some $m\in\mathbb{N}$. 
    \item All $Q^\alpha$'s are negative semidefinite.
\end{enumerate}
Then the $M$-step $N$-stage mixed scheme (\ref{eq:scheme}) is unconditionally energy stable. 

\end{claim}

\begin{proof}
Recall that 
\begin{equation}
    \mathcal{E}(u_{n+1})\leq \mathcal{E}(u_n)+\frac{1}{2k}\sum_{i,j}w_{i,j}d^2(v_i,v_j).
\end{equation}
From Assumption 1, for each $\alpha$ there exists $\{p_j^\alpha\}_{j=-M+1}^N\subset \mathbb{R}^m$ such that $|p_i^\alpha-p_j^\alpha|=d(v_i,v_j)$ for all $(i,j)\in E^\alpha$. It yields
\begin{align}
    \sum_{i,j} w_{i,j}d^2(v_i,v_j) 
    &= \sum_{i,j} \sum_\alpha \theta_{i,j}^\alpha |p_i^\alpha - p_j^\alpha|^2 \nonumber\\
    &=\sum_k\sum_{i,j}\sum_\alpha \theta_{i,j}^\alpha ((p_i^\alpha)_k-(p_j^\alpha)_k)^2\\
    &=\sum_k \sum_\alpha \begin{pmatrix}
     (p_{-M+1}^\alpha)_k\\ \vdots \\ (p_N^\alpha)_k
    \end{pmatrix}^T Q^\alpha \begin{pmatrix}
     (p_{-M+1}^\alpha)_k\\ \vdots \\ (p_N^\alpha)_k
    \end{pmatrix},
    \label{eq:claim1Pf}
\end{align}
where $(p_i^\alpha)_k$ denotes the $k$-th element of $p_j^\alpha \in \mathbb{R}^m$. Since all $Q^\alpha$'s are negative semidefinite, $\mathcal{E}(u_{n+1})\leq \mathcal{E}(u_n)$. 
\end{proof}
Given a graph, there are many ways to partition it into embeddable subgraphs. As an example, we will consider specific types of schemes and suggest possible decomposition which allows us to find concrete instances that satisfy stability and high order consistency concurrently.

\begin{remark}
If Assumption 2 of Claim \ref{claim:stability} is strengthened to require the matrices $Q^\alpha$ to be negative definite, the proof in fact shows that 
    \begin{equation}
        \mathcal{E}(u_{n+1})+\frac{L}{2k}d^2(u_{n+1},u_n)\leq \mathcal{E}(u_n),
    \end{equation}
    for some $L>0$. By induction on $n$, we obtain
    \begin{equation}
        \mathcal{E}(u_{n+1})+\frac{L}{2k}\sum_j d^2(u_{j+1},u_j) \leq \mathcal{E}(u_0).
        \label{eq:strongerStability}
    \end{equation}
    Since the energy $\mathcal{E}$ is bounded from below, it follows that 
    \begin{equation}
        \sum_j d^2(u_{j+1},u_j)<Ck,
        \label{eq:sumd^2}
    \end{equation}
    for some $C>0$. Inequality (\ref{eq:sumd^2}) is a form of equicontinuity in time. Together with boundedness of the energy, this is typically an important step towards a proof of convergence.
    
\end{remark}

We can show the conditions which guarantee energy boundedness of the scheme by slightly modifying $w_{i,j}$ in Claim \ref{claim:stability}. We record this as a relaxed version of stability in the next claim. 

\begin{claim}
\label{claim:boundedness}
Let $(X,d)$ be a metric space and consider the  $M$-step $N$-stage mixed scheme (\ref{eq:scheme}). Recall that $u_j$ denotes the discrete solution of (\ref{eq:gradientFlow}) at the $j$-th time step, and assume $d^2(u_1,u_0)\leq C_0k$ for some constant $C_0>0$. For $0\leq L_1<L_2$, define $\tilde{w}_{i,j}$ as:
\begin{equation}
    \tilde{w}_{i,j}=\begin{dcases}
       w_{i,j}-L_1, \quad &\text{for } i=0, \, j=-1\\
        w_{i,j}+L_2, \quad &\text{for } i=N, \, j=0\\
        w_{i,j}, \quad &\text{elsewhere}
    \end{dcases}
    \label{eq:weightstilde}
\end{equation}
for $i=0,\cdots,N$ and $j=-M+1,\cdots,N$. Replace $w$ by $\tilde{w}$ in Claim 1 and suppose $\{S^\alpha\}$ and $\{Q^\alpha\}$ corresponding to $\tilde{w}$ satisfy Assumption 1 \& 2. Then the energy is uniformly bounded over time.
\end{claim}

\begin{proof}
As in Claim \ref{claim:stability}, we choose $\theta_{i,j}^\alpha \in \mathbb{R}$ such that $\sum_\alpha \theta_{i,j}^\alpha = \tilde{w}_{i,j}$ for each $\alpha$ and $(i,j)\in E^\alpha$. By definition of $\tilde{w}$,
\begin{align}
     &\sum_{i,j}\tilde{w}_{i,j}d^2(v_i,v_j) \\
     &= \sum_{i,j}w_{i,j}d^2(v_i,v_j) - L_1d^2(v_0,v_{-1})+L_2d^2(v_N,v_0) \leq 0.
\end{align}
The last inequality holds from Assumption 1 \& 2. By combining the inequality with (\ref{eq:stabilityIneq}), we obtain
\begin{equation}
     \mathcal{E}(u_{n+1})+\frac{L_2}{2k}d^2(u_{n+1},u_n)\leq \mathcal{E}(u_n)+\frac{L_1}{2k}d^2(u_n,u_{n-1}).
\end{equation}
Since $L_1<L_2$, by induction, 
\begin{equation}
    \mathcal{E}(u_{n+1})+\frac{L_2}{2k}d^2(u_{n+1},u_n)\leq \cdots \leq \mathcal{E}(u_1)+\frac{L_2}{2k}d^2(u_1,u_0).
\end{equation}
It yields
\begin{equation}
    \mathcal{E}(u_{n+1})\leq \mathcal{E}(u_1)+\tilde{C}
\end{equation}
for some constant $\tilde{C}\in\mathbb{R}$. 

\end{proof}

We now consider two examples of embeddable graphs that will be used, in conjuction with Claims 1 \& 2, to seek multi-stage multi-step numerical schemes of the form (\ref{eq:scheme}) that have provable unconditional stability properties in any metric space.

\textbf{Example. Multi-stage schemes}

We consider the single-step $N$-stage scheme of the form (\ref{eq:scheme}). Note that $\gamma_{i,j}$'s appear on the scheme for $i=1,\cdots,N$ and $j=0,\cdots,i-1$. The values of  $w_{i,j}$ are defined for $i,j=0,\cdots,N$ by formula (\ref{eq:weights}). Based on the values of $\gamma_{i,j}$, the associated graph should be either a complete undirected graph on $(N+1)$ vertices or a subgraph of it. Let $G$ be a complete undirected graph on $V=\{0,1,\cdots,N\}$; see Figure \ref{fig:NstageGraph}. Then $G$ can be decomposed into $(N-1)$ subgraphs as follows: $G=\cup_{\alpha=0}^{N-2}S^\alpha$ where $S^\alpha=(V,E^\alpha)$ and 
\begin{align}
    E^\alpha = &\{(i,j): i=\alpha,\alpha+1, j=i+1,\cdots,N\} \nonumber\\
    &\cup \{(i,j):i=j+1,\cdots,N,j =\alpha,\alpha+1\}.
    \label{eq:NstageSubgraph}
\end{align}
Subgraphs $S^0$, $S^1$ and $S^2$ are illustrated in Figure \ref{fig:NstageSubgraph}. Note that the subgraphs defined as above are embeddable into $\mathbb{R}^2$. Suppose that all distances are specified. For two distinguished vertices $i$ and $j$, there exist $p_i$ and $p_j$ in $\mathbb{R}^2$ such that $|p_i-p_j|=d_{i,j}$. For every triplet $(d_{i,j},d_{i,k},d_{j,k})$, we can choose $p_k\in\mathbb{R}^2$ such that $|p_i-p_k|=d_{i,k}$ and $|p_j-p_k|=d_{j,k}$ since there exists a triangle on $\mathbb{R}^2$ whose sides are equal to each element of the triplet. All of triangles in $S^\alpha$ have a common edge, therefore the repetitive use of this method shows the embeddability. 

\begin{figure}[h!]
\centering
\begin{tikzpicture}
    \node[draw=none, minimum size = 4cm, regular polygon, regular polygon sides = 8] (a) {};
    \foreach \x in {1,2,...,5} 
        \fill (a.corner \x) circle [radius=2pt];
    \foreach \x in {7,8} 
    \fill (a.corner \x) circle [radius=2pt];
    \draw (a.corner 1) node[above] {$0$};
    \draw (a.corner 2) node[above] {$1$};
    \draw (a.corner 3) node[left] {$2$};
    \draw (a.corner 4) node[below] {$3$};
    \draw (a.corner 5) node[below] {$4$};
    \draw (a.corner 7) node[below] {${N-1}$};
    \draw (a.corner 8) node[right] {${N}$};
    \draw (a.corner 1)--(a.corner 2);
    \draw (a.corner 1)--(a.corner 3);
    \draw (a.corner 1)--(a.corner 4);
    \draw (a.corner 1)--(a.corner 5);
    \draw (a.corner 1)--(a.corner 7);
    \draw (a.corner 1)--(a.corner 8);
    \draw (a.corner 2)--(a.corner 3);
    \draw (a.corner 2)--(a.corner 4);
    \draw (a.corner 2)--(a.corner 5);
    \draw (a.corner 2)--(a.corner 7);
    \draw (a.corner 2)--(a.corner 8);
    \draw (a.corner 3)--(a.corner 4);
    \draw (a.corner 3)--(a.corner 5);
    \draw (a.corner 3)--(a.corner 7);
    \draw (a.corner 3)--(a.corner 8);
    \draw (a.corner 4)--(a.corner 5);
    \draw (a.corner 4)--(a.corner 8);
    \draw (a.corner 7)--(a.corner 8);
    \draw (a.corner 5) node[above=0.5cm, right=1cm] {$\iddots$};
\end{tikzpicture}
\caption{The graph $G$, a complete undirected graph on a set of vertices $\{0,1,\cdots,N\}$, associated with the single-step $N$-stage scheme of the form (\ref{eq:scheme}).} 
\label{fig:NstageGraph}
\end{figure}
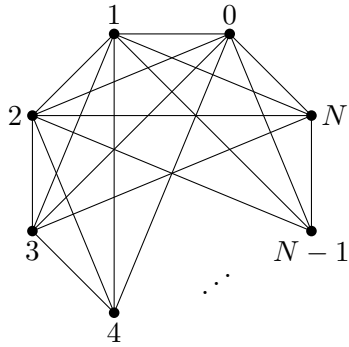

\begin{figure}[h!]
\centering
\begin{subfigure}[b]{0.3\linewidth}
    \centering
    \begin{tikzpicture}
    \node[draw=none, minimum size = 3cm, regular polygon, regular polygon sides = 8] (a) {};
    \foreach \x in {1,2,...,5} 
        \fill (a.corner \x) circle [radius=2pt];
    \foreach \x in {7,8} 
    \fill (a.corner \x) circle [radius=2pt];
    \draw (a.corner 1) node[above] {$0$};
    \draw (a.corner 2) node[above] {$1$};
    \draw (a.corner 3) node[left] {$2$};
    \draw (a.corner 4) node[below] {$3$};
    \draw (a.corner 5) node[below] {$4$};
    \draw (a.corner 7) node[below] {${N-1}$};
    \draw (a.corner 8) node[right] {${N}$};
    \draw (a.corner 1)--(a.corner 2);
    \draw (a.corner 1)--(a.corner 3);
    \draw (a.corner 1)--(a.corner 4);
    \draw (a.corner 1)--(a.corner 5);
    \draw (a.corner 1)--(a.corner 7);
    \draw (a.corner 1)--(a.corner 8);
    \draw (a.corner 2)--(a.corner 3);
    \draw (a.corner 2)--(a.corner 4);
    \draw (a.corner 2)--(a.corner 5);
    \draw (a.corner 2)--(a.corner 7);
    \draw (a.corner 2)--(a.corner 8);
    \draw (a.corner 5) node[above=0.4cm, right=0.5cm] {$\iddots$};
\end{tikzpicture}
\caption{$S^0$}
\end{subfigure}
\begin{subfigure}[b]{0.3\linewidth}
    \centering
    \begin{tikzpicture}
    \node[draw=none, minimum size = 3cm, regular polygon, regular polygon sides = 8] (a) {};
    \foreach \x in {1,2,...,5} 
        \fill (a.corner \x) circle [radius=2pt];
    \foreach \x in {7,8} 
    \fill (a.corner \x) circle [radius=2pt];
    \draw (a.corner 1) node[above] {$0$};
    \draw (a.corner 2) node[above] {$1$};
    \draw (a.corner 3) node[left] {$2$};
    \draw (a.corner 4) node[below] {$3$};
    \draw (a.corner 5) node[below] {$4$};
    \draw (a.corner 7) node[below] {${N-1}$};
    \draw (a.corner 8) node[right] {${N}$};
    \draw (a.corner 2)--(a.corner 3);
    \draw (a.corner 2)--(a.corner 4);
    \draw (a.corner 2)--(a.corner 5);
    \draw (a.corner 2)--(a.corner 7);
    \draw (a.corner 2)--(a.corner 8);
    \draw (a.corner 3)--(a.corner 4);
    \draw (a.corner 3)--(a.corner 5);
    \draw (a.corner 3)--(a.corner 7);
    \draw (a.corner 3)--(a.corner 8);
    \draw (a.corner 5) node[above=0.4cm, right=0.5cm] {$\iddots$};
\end{tikzpicture}
\caption{$S^1$}
\end{subfigure}
\begin{subfigure}[b]{0.3\linewidth}
    \centering
    \begin{tikzpicture}
    \node[draw=none, minimum size = 3cm, regular polygon, regular polygon sides = 8] (a) {};
    \foreach \x in {1,2,...,5} 
        \fill (a.corner \x) circle [radius=2pt];
    \foreach \x in {7,8} 
    \fill (a.corner \x) circle [radius=2pt];
    \draw (a.corner 1) node[above] {$0$};
    \draw (a.corner 2) node[above] {$1$};
    \draw (a.corner 3) node[left] {$2$};
    \draw (a.corner 4) node[below] {$3$};
    \draw (a.corner 5) node[below] {$4$};
    \draw (a.corner 7) node[below] {${N-1}$};
    \draw (a.corner 8) node[right] {${N}$};
    \draw (a.corner 3)--(a.corner 4);
    \draw (a.corner 3)--(a.corner 5);
    \draw (a.corner 3)--(a.corner 7);
    \draw (a.corner 3)--(a.corner 8);
    \draw (a.corner 4)--(a.corner 5);
    \draw (a.corner 4)--(a.corner 7);
    \draw (a.corner 4)--(a.corner 8);
    \draw (a.corner 5) node[above=0.4cm, right=0.5cm] {$\iddots$};
\end{tikzpicture}
\caption{$S^2$}
\end{subfigure}
\caption{Subgraphs of $G$ in Figure \ref{fig:NstageGraph} which are embeddable in $\mathbb{R}^2$.}
\label{fig:NstageSubgraph}x
\end{figure}
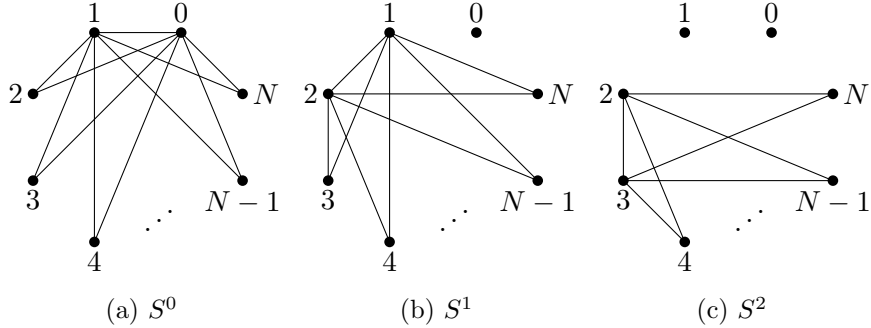

\textbf{Example. Two-step multi-stage diagonal schemes}

\begin{figure}[h!]
\centering
    \begin{tikzpicture}
        \node[draw=none, minimum size = 4cm, regular polygon, regular polygon sides = 8] (a) {};
        \foreach \x in {1,...,8} 
            \fill (a.corner \x) circle [radius=2pt];
        \draw (a.corner 1) node[above] {$0$};
        \draw (a.corner 2) node[left] {$-1$};
        \draw (a.corner 3) node[left] {$1$};
        \draw (a.corner 4) node[left] {$2$};
        \draw (a.corner 5) node[below] {$3$};
        \draw (a.corner 6) node[below] {$4$};
        \draw (a.corner 7) node[right] {$N-1$};
        \draw (a.corner 8) node[right] {$N$};
        \draw (a.corner 1)--(a.corner 2);
        \draw (a.corner 1)--(a.corner 3);
        \draw (a.corner 1)--(a.corner 4);
        \draw (a.corner 1)--(a.corner 5);
        \draw (a.corner 1)--(a.corner 6);
        \draw (a.corner 1)--(a.corner 7);
        \draw (a.corner 1)--(a.corner 8);
        \draw (a.corner 2)--(a.corner 3);
        \draw (a.corner 2)--(a.corner 4);
        \draw (a.corner 2)--(a.corner 5);
        \draw (a.corner 2)--(a.corner 6);
        \draw (a.corner 2)--(a.corner 7);
        \draw (a.corner 2)--(a.corner 8);
        \draw (a.corner 3)--(a.corner 4);
        \draw (a.corner 4)--(a.corner 5);
        \draw (a.corner 5)--(a.corner 6);
        \draw (a.corner 7)--(a.corner 8);
        \draw (a.corner 6) node[above=0.6cm, right=0.4cm] {$\iddots$};
    \end{tikzpicture}
    \caption{The graph $G$ associated with the $2$-step $N$-stage diagonal scheme of the form (\ref{eq:scheme}).}
    \label{fig:2stepNstageGraph}
\end{figure}
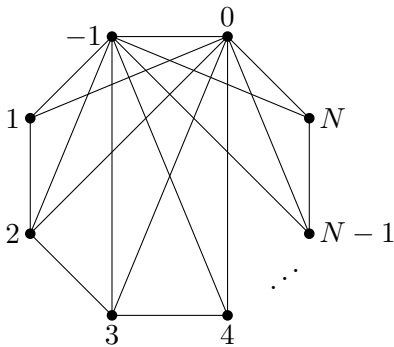

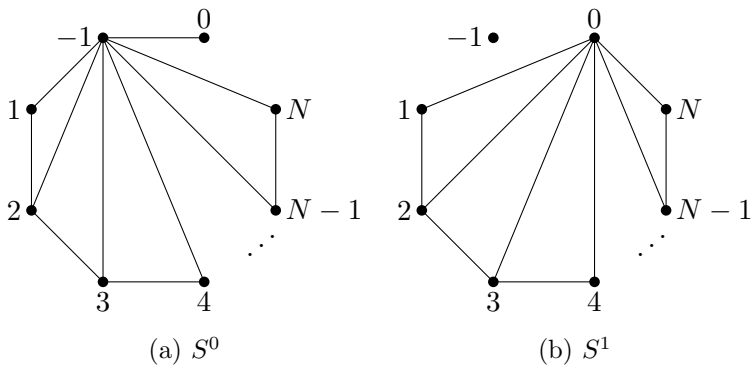
\begin{figure}
    \centering
    \begin{subfigure}[b]{0.4\linewidth}
    \centering
    \begin{tikzpicture}
        \node[draw=none, minimum size = 3.5cm, regular polygon, regular polygon sides = 8] (a) {};
        \foreach \x in {1,...,8} 
            \fill (a.corner \x) circle [radius=2pt];
        \draw (a.corner 1) node[above] {$0$};
        \draw (a.corner 2) node[left] {$-1$};
        \draw (a.corner 3) node[left] {$1$};
        \draw (a.corner 4) node[left] {$2$};
        \draw (a.corner 5) node[below] {$3$};
        \draw (a.corner 6) node[below] {$4$};
        \draw (a.corner 7) node[right] {$N-1$};
        \draw (a.corner 8) node[right] {$N$};
        \draw (a.corner 1)--(a.corner 2);
        \draw (a.corner 2)--(a.corner 3);
        \draw (a.corner 2)--(a.corner 4);
        \draw (a.corner 2)--(a.corner 5);
        \draw (a.corner 2)--(a.corner 6);
        \draw (a.corner 2)--(a.corner 7);
        \draw (a.corner 2)--(a.corner 8);
        \draw (a.corner 3)--(a.corner 4);
        \draw (a.corner 4)--(a.corner 5);
        \draw (a.corner 5)--(a.corner 6);
        \draw (a.corner 7)--(a.corner 8);
        \draw (a.corner 6) node[above=0.6cm, right=0.4cm] {$\iddots$};
    \end{tikzpicture}
    \caption{$S^0$}
    \end{subfigure}
    \begin{subfigure}[b]{0.4\linewidth}
    \centering
    \begin{tikzpicture}
        \node[draw=none, minimum size = 3.5cm, regular polygon, regular polygon sides = 8] (a) {};
        \foreach \x in {1,...,8} 
            \fill (a.corner \x) circle [radius=2pt];
        \draw (a.corner 1) node[above] {$0$};
        \draw (a.corner 2) node[left] {$-1$};
        \draw (a.corner 3) node[left] {$1$};
        \draw (a.corner 4) node[left] {$2$};
        \draw (a.corner 5) node[below] {$3$};
        \draw (a.corner 6) node[below] {$4$};
        \draw (a.corner 7) node[right] {$N-1$};
        \draw (a.corner 8) node[right] {$N$};
        \draw (a.corner 1)--(a.corner 3);
        \draw (a.corner 1)--(a.corner 4);
        \draw (a.corner 1)--(a.corner 5);
        \draw (a.corner 1)--(a.corner 6);
        \draw (a.corner 1)--(a.corner 7);
        \draw (a.corner 1)--(a.corner 8);
        \draw (a.corner 3)--(a.corner 4);
        \draw (a.corner 4)--(a.corner 5);
        \draw (a.corner 5)--(a.corner 6);
        \draw (a.corner 7)--(a.corner 8);
        \draw (a.corner 6) node[above=0.6cm, right=0.4cm] {$\iddots$};
    \end{tikzpicture}
    \caption{$S^1$}
    \end{subfigure}
    \caption{Subgraphs of $G$ in Figure \ref{fig:2stepNstageGraph} which are embeddable in $\mathbb{R}^2$.}
    \label{fig:2stepNstageSubgraph}
\end{figure}

We now focus on a special case of multi-step, multi-stage (mixed) variational schemes of the form (\ref{eq:scheme}), where in each optimization problem the movement limiter terms depend only on the most recent intermediate stage (but may depend on multiple previous steps). We call it a multi-stage diagonal scheme. As an example, we consider the $2$-step $N$-stage diagonal scheme of the form (\ref{eq:scheme}). Since $\gamma_{i,j}=0$ for $j=0,\cdots,i-2$, due to the definition of a multi-step diagonal scheme, the graph corresponding to this scheme connects the  $i$-th vertex to the $(-1)$-th, $0$-th, and $(i-1)$-th vertices respectively. (Or, it should be a subgraph of it depending on the values of $\gamma_{i,j}$.) 
The graph $G$ defined on a set of vertices $V=\{-1,0,\cdots,N\}$ associated with this scheme is represented in Figure \ref{fig:2stepNstageGraph}. 
$G$ can be partitioned into two embeddable subgraphs $S^0$ and $S^1$ where $S^\alpha=(V,E^\alpha)$ and
\begin{align}
    E^\alpha =  &\{(\alpha-1,j):j=\alpha,\cdots,N\}\cup \{(j,\alpha-1):j=\alpha,\cdots,N\} \nonumber\\
    &\cup \{(j+1,j):j=1,\cdots,N\}\cup \{(j,j+1):j=1,\cdots,N\}.
    \label{eq:2stepNstageSubgraph}
\end{align}
for $\alpha=0,1$; see Figure \ref{fig:2stepNstageSubgraph}. 
 
\begin{remark}
    For the coefficients $\gamma_{i,j}$'s that appear in scheme (\ref{eq:scheme}), we have exhibited stability conditions that are independent of the time step size $k$. However, in general, well-posedness (existence and uniqueness of a global minimizer) of the variational problems (\ref{eq:scheme}) that need to be solved at each step of the scheme may require conditions to be placed on the time step size. In addition, there is a much more serious issue of {\em how} to solve these optimization problems reliably and efficiently. All these important issues need to be addressed on a problem by problem basis. 
\end{remark}

\section{Consistency}
\label{sec:consistency}

The Jordan-Kinderlehrer-Otto scheme \cite{jko} generates a numerical solution for gradient flow (\ref{eq:gradientFlow}) by solving the variational problem: 
\begin{equation}
    \min_{\xi\in X} \left\{\mathcal{E}(\xi)+\frac{1}{2k}D(\xi,u_n)\right\}
    \label{eq:JKO}
\end{equation}
where $D=d^2$. This scheme is known to be first order accurate in time. We now make an ansatz for the next step $u_{n+1}$ as $u_{n+1}=u_n+qk+O(k^2)$ for some $q\in X$, and plug it into the Euler-Lagrange equation of the optimization problem (\ref{eq:JKO}). From the first order term in $k$, we obtain the first order time derivative of the exact solution as follows:
\begin{equation}
    u_t =-2(\nabla^2 D(u,u))^{-1}\nabla \mathcal{E}(u).
    \label{eq:u_t}
\end{equation}

As our scheme (\ref{eq:scheme}) allows us to use multi-step and multi-stage, it is possible to generate a numerical solution which matches the exact solution up to higher order in $k$. We will present the conditions on $\gamma_{i,j}$'s which guarantee high order accuracy of the scheme in the next claim, and add the proof in Section $\ref{sec:Appendix}$ (Appendix). 

\begin{claim}
\label{claim:consistency}
Given a set of coefficients \{$\gamma_{i,j}\in \mathbb{R}$: $i=1,\cdots,N$, $j=-M+1,\cdots, i-1$\}, define auxiliary variables as follows: 
\begin{equation}
\begin{dcases}
     a_i &= \frac{1}{\sum_j \gamma_{i,j}}\left(1+\sum_j \gamma_{i,j}a_j\right)\\
     b_i &= \frac{1}{\sum_j \gamma_{i,j}} \left(a_i+\sum_j \gamma_{i,j}b_j\right)\\
     c_i &= \frac{1}{\sum_{j}\gamma_{i,j}}\left(b_i + \sum_{j}\gamma_{i,j}c_j\right)\\
     d_i &=\frac{1}{\sum_{j}\gamma_{i,j}}\left(\frac{a_i^2}{2} + \sum_{j}\gamma_{i,j}d_j\right).
     \label{eq:consistencyVars}
\end{dcases}
\end{equation}
The $M$-step $N$-stage scheme (\ref{eq:scheme}) is second order accurate when
\begin{equation}
    a_N=1, \quad b_N=\frac{1}{2},
    \label{eq:2ndConsistency}
\end{equation}
and third  order accurate when
\begin{equation}
    a_N=1, \quad b_N=\frac{1}{2}, \quad c_N=\frac{1}{6},\quad d_N=\frac{1}{6}.
    \label{eq:3rdConsistency}
    \end{equation}
\end{claim}
\begin{proof}
See Appendix. 
\end{proof}

We now exhibit concrete examples of high order accurate schemes which are either unconditionally energy stable or energy bounded. 

\textbf{Example. Energy stable second order scheme}

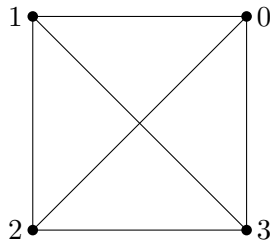
\begin{figure}[h!]
\centering
    \begin{tikzpicture}
        \node[draw=none, minimum size = 4cm, regular polygon, regular polygon sides = 4] (a) {};
        \foreach \x in {1,...,4} 
            \fill (a.corner \x) circle [radius=2pt];
        \draw (a.corner 1) node[right] {$0$};
        \draw (a.corner 2) node[left] {$1$};
        \draw (a.corner 3) node[left] {$2$};
        \draw (a.corner 4) node[right] {$3$};
        \draw (a.corner 1)--(a.corner 2);
        \draw (a.corner 1)--(a.corner 3);
        \draw (a.corner 1)--(a.corner 4);
        \draw (a.corner 2)--(a.corner 3);
        \draw (a.corner 2)--(a.corner 4);
        \draw (a.corner 3)--(a.corner 4);
    \end{tikzpicture}
    \caption{A complete undirected graph on $4$ vertices $\{0,1,2,3\}$ associated with the  single-step $3$-stage scheme (\ref{eq:2ndScheme}).}
    \label{fig:3stageGraph}
\end{figure}

As an example of the second order scheme, we consider the single-step $3$-stage scheme with $\gamma_{i,j}$'s given by: 
\begin{equation}
   \begin{pmatrix}
    \gamma_{1,0} & 0 & 0\\ 
     \gamma_{2,0} & \gamma_{2,1} & 0\\ 
      \gamma_{3,0} & \gamma_{3,1} & \gamma_{3,2}
   \end{pmatrix}=  \begin{pmatrix}
    4 & 0 & 0 \\
    -1 & 5 & 0 \\
    -2 & -1.6 & 9.6
    \end{pmatrix}.
    \label{eq:2ndScheme}
\end{equation}
We can easily see that the corresponding auxiliary variables $a_3$ and $b_3$ in formula (\ref{eq:consistencyVars}) satisfy the second order consistency conditions (\ref{eq:2ndConsistency}). To prove the stability of the scheme, we find the values of $w_{i,j}$'s using formula (\ref{eq:weights}):
\begin{equation}
    \begin{pmatrix}
    w_{1,0} & 0 & 0\\ 
     w_{2,0} & w_{2,1} & 0\\ 
      w_{3,0} & w_{3,1} & w_{3,2}
   \end{pmatrix}=  \begin{pmatrix}
    -5 & 0 & 0 \\
    -1 & -6.6 & 0 \\
    2 & 1.6 & -9.6
    \end{pmatrix},
    \label{eq:2ndSchemeWeights}
\end{equation}
and all values of $w_{i,j}$'s not shown in the matrix (\ref{eq:2ndSchemeWeights}) are zero. According to the values of $w_{i,j}$'s, the graph associated with this scheme is a complete undirected graph on $4$ vertices $V=\{0,1,2,3\}$; see Figure \ref{fig:3stageGraph}.

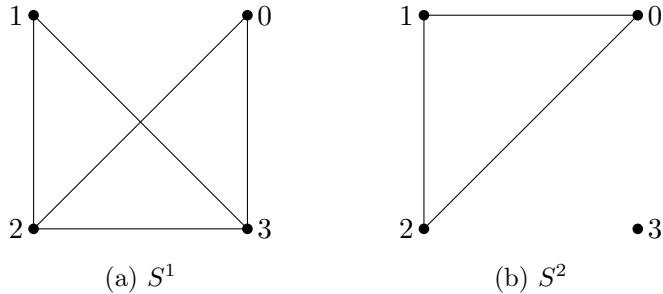
\begin{figure}[h!]
\centering
\begin{subfigure}[b]{0.4\linewidth}
    \centering
    \begin{tikzpicture}
        \node[draw=none, minimum size = 4cm, regular polygon, regular polygon sides = 4] (a) {};
        \foreach \x in {1,...,4} 
            \fill (a.corner \x) circle [radius=2pt];
        \draw (a.corner 1) node[right] {$0$};
        \draw (a.corner 2) node[left] {$1$};
        \draw (a.corner 3) node[left] {$2$};
        \draw (a.corner 4) node[right] {$3$};
        \draw (a.corner 1)--(a.corner 3);
        \draw (a.corner 1)--(a.corner 4);
        \draw (a.corner 2)--(a.corner 3);
        \draw (a.corner 2)--(a.corner 4);
        \draw (a.corner 3)--(a.corner 4);
    \end{tikzpicture}
    \caption{$S^1$}
\end{subfigure}
\begin{subfigure}[b]{0.4\linewidth}
    \centering
    \begin{tikzpicture}
        \node[draw=none, minimum size = 4cm, regular polygon, regular polygon sides = 4] (a) {};
        \foreach \x in {1,...,4} 
            \fill (a.corner \x) circle [radius=2pt];
        \draw (a.corner 1) node[right] {$0$};
        \draw (a.corner 2) node[left] {$1$};
        \draw (a.corner 3) node[left] {$2$};
        \draw (a.corner 4) node[right] {$3$};
       \draw (a.corner 1)--(a.corner 2);
       \draw (a.corner 1)--(a.corner 3);
        \draw (a.corner 2)--(a.corner 3);
    \end{tikzpicture}
    \caption{$S^2$}
\end{subfigure}
\caption{Subgraphs of a graph in Figure \ref{fig:3stageGraph} which are embeddable into $\mathbb{R}^2$.}
\label{fig:3stageSubgraph}
\end{figure}

We decompose the graph into embeddable subgraphs $S^1=(V,E^1)$ and $S^2=(V,E^2)$ where 
\begin{equation}
    E^1 = \{(0,2),(2,0),(0,3),(3,0),(1,2),(2,1),(1,3),(3,1),(2,3),(3,2)\}
\end{equation}
and
\begin{equation}
    E^2 = \{(0,1),(1,0),(0,2),(2,0),(1,2),(2,1)\};
\end{equation}
see Figure \ref{fig:3stageSubgraph}. For each subgraph, choose $\theta^1_{i,j}$'s and $\theta^2_{i,j}$'s as follows:
\begin{align}
    \begin{pmatrix}
    0 & 0 & 0\\ 
    \theta^1_{2,0} & \theta^1_{2,1} & 0\\ 
    \theta^1_{3,0} & \theta^1_{3,1} & \theta^1_{3,2}
   \end{pmatrix} &=  \begin{pmatrix}
    0 & 0 & 0 \\
    -3.2 & -2.67 & 0 \\
    2 & 1.6 & -9.6
    \end{pmatrix},\\
    \begin{pmatrix}
    \theta^2_{1,0} & 0 & 0\\ 
    \theta^2_{2,0} & \theta^2_{2,1} & 0\\ 
    0 & 0 & 0
   \end{pmatrix} &=  \begin{pmatrix}
    -5 & 0 & 0 \\
    2.2 & -3.93 & 0 \\
    0 & 0 & 0
    \end{pmatrix},
\end{align}
and all $\theta^1_{i,j}$'s and $\theta^2_{i,j}$'s not shown above are zero. Then 
\begin{equation}
    \theta^1_{i,j}+\theta^2_{i,j}=w_{i,j}
\end{equation}
for all $i,j$, and we can check that eigenvalues of matrices $Q^1$ and $Q^2$ defined as 
\begin{equation}
    Q^\alpha = \sum_{(i,j)\in E^\alpha} \theta_{i,j}^\alpha (e_{i+1}-e_{j+1})\otimes (e_{i+1}-e_{j+1})
\end{equation}
for $\alpha=1,2$ are all non-positve. Therefore, by Claim \ref{claim:stability}, the single-step $3$-stage scheme (\ref{eq:2ndScheme}) is energy stable.

\textbf{Example. Energy bounded third  order scheme}

To obtain higher order accuracy, we are concerned with multi-step multi-stage mixed schemes. Let us consider the $2$-step $7$-stage diagonal scheme with $\gamma_{i,j}$'s given by:
\begin{equation}
    \begin{pmatrix}
     \gamma_{1,-1} & \gamma_{1,0} & 0 \\
    \gamma_{2,-1} & \gamma_{2,0} & \gamma_{2,1} \\
    \vdots\\
    \gamma_{6,-1} & \gamma_{6,0} & \gamma_{6,5} \\
    \gamma_{7,-1} & \gamma_{7,0} & \gamma_{7,6} \\
    \end{pmatrix}
    \approx 
\begin{pmatrix}
 0.20  & 12.96    &     0 \\
  -0.67 &   0.64 &  12.45\\
  -0.01  & -0.76  & 13.27\\
    0.26  & -1.42 &   8.97\\
    0.05  & -0.62 &   6.90\\
    0.01  & -1.49  &  8.31\\
    0.19 &  -0.59  & 11.25
\end{pmatrix},
\label{eq:3rdScheme}
\end{equation}
and $\gamma_{i,j}=0$ for $i=3,\cdots,7$ and $j=1,\cdots,i-2$. The exact values of $\gamma_{i,j}$'s are recorded in Section \ref{sec:Appendix} (Appendix). 

This scheme is third  order accurate as the associated auxiliary variables $a_7,b_7,c_7$ and $d_7$ in formula (\ref{eq:consistencyVars}) satisfy the third order consistency conditions  (\ref{eq:3rdConsistency}). We now show this choice of $\gamma_{i,j}$'s  ensures the uniform boundedness of the energy. Choose $L_1 =0.2$ and $L_2=0.3$, and find the values of $\tilde{w}_{i,j}$'s in formula (\ref{eq:weightstilde}):
\begin{equation}
   \begin{pmatrix}
    \tilde{w}_{1,-1} & \tilde{w}_{1,0} & 0 \\
    \tilde{w}_{2,-1} & \tilde{w}_{2,0} & \tilde{w}_{2,1} \\
    \tilde{w}_{3,-1} & \tilde{w}_{3,0} & \tilde{w}_{3,1} \\
    \vdots&\vdots \\
    \tilde{w}_{6,-1} & \tilde{w}_{6,0} & \tilde{w}_{6,5} \\
    \tilde{w}_{7,-1} & \tilde{w}_{7,0} & \tilde{w}_{7,6} \\
    \end{pmatrix}
    \approx 
\begin{pmatrix}
 -0.87  & -12.32    & 0      \\
  0.66 &   -1.40 &  -12.45\\
  0.27  & -0.66  & -13.27\\
    -0.21  & 0.80 &  -8.97\\
    -0.04  & -0.87 &  -6.90\\
    0.17  & 0.90  &  -8.31\\
    -0.19 &  0.89  & -11.25
\end{pmatrix},
\end{equation}
and all $\tilde{w}_{i,j}$'s not shown above are zero. We partition the graph associated with this scheme into two embeddable subgraphs $S^1$ and $S^2$ as in Figure \ref{fig:2step7stageSubgraph}, and choose $\theta_{i,j}^1$ and $\theta_{i,j}^2$'s as follows: 
\begin{equation}
    \begin{pmatrix}
    \theta^1_{1,-1} & 0 \\
    \theta^1_{2,-1} & \theta^1_{2,1} \\
    \vdots\\
    \theta^1_{6,-1} & \theta^1_{6,5} \\
    \theta^1_{7,-1} & \theta^1_{7,6} \\
    \end{pmatrix}
    \approx 
\begin{pmatrix}
 -0.87     &     0 \\
   0.66  &  -4.82\\
  0.27   & -2.46\\
    -0.21  &   -1.52\\
    -0.04  &   -0.10\\
    0.17   &  -0.12\\
    -0.19   & -0.94
\end{pmatrix}, \quad 
 \begin{pmatrix}
    \theta^2_{1,0} & 0 \\
    \theta^2_{2,0} & \theta^2_{2,1} \\
    \vdots\\
    \theta^2_{6,0} & \theta^2_{6,5} \\
    \theta^2_{7,0} & \theta^2_{7,6} \\
    \end{pmatrix}
    \approx 
\begin{pmatrix}
 -12.32    &     0 \\
   -1.40 &  -7.63\\
  -0.66  & -10.81\\
    0.80 &   -7.45\\
    -0.87 &   -6.80\\
    0.90  &  -8.18\\
    0.89  & -10.31
\end{pmatrix},
\end{equation}
and all $\theta_{i,j}^1$'s and $\theta_{i,j}^2$'s not shown above are zero. This choice of $\theta_{i,j}^1$'s and $\theta_{i,j}^2$'s fulfills Assumption $2$ in Claim \ref{claim:stability}, and hence, the $2$-step $7$-stage diagonal scheme (\ref{eq:3rdScheme}) ensures the energy boundedness over time by Claim \ref{claim:boundedness}.

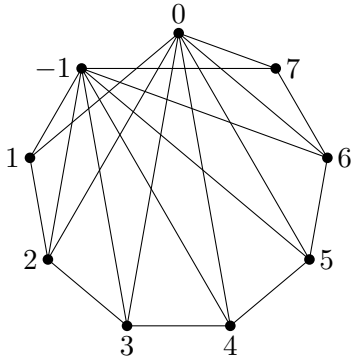
\begin{figure}[h!]
\centering
    \begin{tikzpicture}
        \node[draw=none, minimum size = 4cm, regular polygon, regular polygon sides = 9] (a) {};
        \foreach \x in {1,...,9} 
            \fill (a.corner \x) circle [radius=2pt];
        \draw (a.corner 1) node[above] {$0$};
        \draw (a.corner 2) node[left] {$-1$};
        \draw (a.corner 3) node[left] {$1$};
        \draw (a.corner 4) node[left] {$2$};
        \draw (a.corner 5) node[below] {$3$};
        \draw (a.corner 6) node[below] {$4$};
        \draw (a.corner 7) node[right] {$5$};
        \draw (a.corner 8) node[right] {$6$};
        \draw (a.corner 9) node[right] {$7$};
        \draw (a.corner 1)--(a.corner 3);
        \draw (a.corner 1)--(a.corner 4);
        \draw (a.corner 1)--(a.corner 5);
        \draw (a.corner 1)--(a.corner 6);
        \draw (a.corner 1)--(a.corner 7);
        \draw (a.corner 1)--(a.corner 8);
        \draw (a.corner 1)--(a.corner 9);
        \draw (a.corner 2)--(a.corner 3);
        \draw (a.corner 2)--(a.corner 4);
        \draw (a.corner 2)--(a.corner 5);
        \draw (a.corner 2)--(a.corner 6);
        \draw (a.corner 2)--(a.corner 7);
        \draw (a.corner 2)--(a.corner 8);
        \draw (a.corner 2)--(a.corner 9);
        \draw (a.corner 3)--(a.corner 4);
        \draw (a.corner 4)--(a.corner 5);
        \draw (a.corner 5)--(a.corner 6);
        \draw (a.corner 6)--(a.corner 7);
        \draw (a.corner 7)--(a.corner 8);
        \draw (a.corner 8)--(a.corner 9);
    \end{tikzpicture}
    \caption{An undirected graph on $8$ vertices $\{-1,0,\cdots,7\}$ associated with the $2$-step $7$-stage diagonal scheme (\ref{eq:3rdScheme}).}
    \label{fig:2step7stageGraph}
\end{figure}

\begin{figure}
    \centering
    \begin{subfigure}[b]{0.4\linewidth}
    \centering
    \begin{tikzpicture}
        \node[draw=none, minimum size = 3.5cm, regular polygon, regular polygon sides = 9] (a) {};
        \foreach \x in {1,...,9} 
            \fill (a.corner \x) circle [radius=2pt];
        \draw (a.corner 1) node[above] {$0$};
        \draw (a.corner 2) node[left] {$-1$};
        \draw (a.corner 3) node[left] {$1$};
        \draw (a.corner 4) node[left] {$2$};
        \draw (a.corner 5) node[below] {$3$};
        \draw (a.corner 6) node[below] {$4$};
        \draw (a.corner 7) node[right] {$5$};
        \draw (a.corner 8) node[right] {$6$};
        \draw (a.corner 9) node[right] {$7$};
        \draw (a.corner 2)--(a.corner 3);
        \draw (a.corner 2)--(a.corner 4);
        \draw (a.corner 2)--(a.corner 5);
        \draw (a.corner 2)--(a.corner 6);
        \draw (a.corner 2)--(a.corner 7);
        \draw (a.corner 2)--(a.corner 8);
        \draw (a.corner 2)--(a.corner 9);
        \draw (a.corner 3)--(a.corner 4);
        \draw (a.corner 4)--(a.corner 5);
        \draw (a.corner 5)--(a.corner 6);
        \draw (a.corner 6)--(a.corner 7);
        \draw (a.corner 7)--(a.corner 8);
        \draw (a.corner 8)--(a.corner 9);
    \end{tikzpicture}
    \caption{$S^1$}
    \end{subfigure}
    \begin{subfigure}[b]{0.4\linewidth}
    \centering
    \begin{tikzpicture}
        \node[draw=none, minimum size = 3.5cm, regular polygon, regular polygon sides = 9] (a) {};
        \foreach \x in {1,...,9} 
            \fill (a.corner \x) circle [radius=2pt];
        \draw (a.corner 1) node[above] {$0$};
        \draw (a.corner 2) node[left] {$-1$};
        \draw (a.corner 3) node[left] {$1$};
        \draw (a.corner 4) node[left] {$2$};
        \draw (a.corner 5) node[below] {$3$};
        \draw (a.corner 6) node[below] {$4$};
        \draw (a.corner 7) node[right] {$5$};
        \draw (a.corner 8) node[right] {$6$};
        \draw (a.corner 9) node[right] {$7$};
        \draw (a.corner 1)--(a.corner 3);
        \draw (a.corner 1)--(a.corner 4);
        \draw (a.corner 1)--(a.corner 5);
        \draw (a.corner 1)--(a.corner 6);
        \draw (a.corner 1)--(a.corner 7);
        \draw (a.corner 1)--(a.corner 8);
        \draw (a.corner 1)--(a.corner 9);
        \draw (a.corner 3)--(a.corner 4);
        \draw (a.corner 4)--(a.corner 5);
        \draw (a.corner 5)--(a.corner 6);
        \draw (a.corner 6)--(a.corner 7);
        \draw (a.corner 7)--(a.corner 8);
        \draw (a.corner 8)--(a.corner 9);
    \end{tikzpicture}
    \caption{$S^2$}
    \end{subfigure}
    \caption{Subgraphs of a graph in Figure \ref{fig:2step7stageGraph} which are embeddable (in $\mathbb{R}^2$).}
    \label{fig:2step7stageSubgraph}
\end{figure}
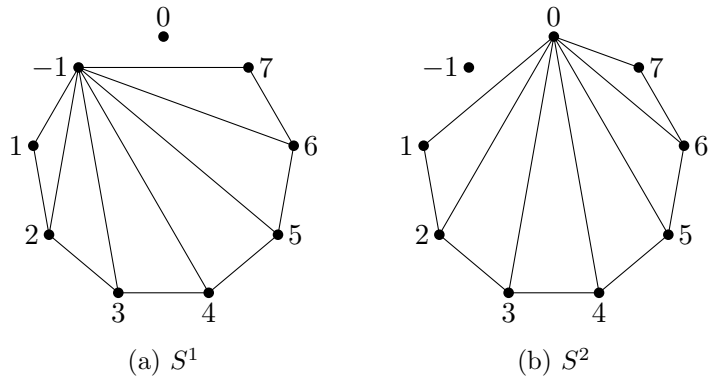

\section{Numerical results}
\label{sec:numResults}

In this section, we apply two schemes  (\ref{eq:2ndScheme}) and (\ref{eq:3rdScheme}) to several gradient flows and demonstrate the accuracy of each scheme. We will consider the  one-dimensional  heat equation, porous medium equation, and nonlinear Fokker-Planck equation which can be seen as gradient flows in the quadratic Wasserstein space. In each of the following examples, the computation domain is $[-1,1]$ and we choose the initial condition $u(x,0)=\frac{1}{2}+\frac{1}{4}\cos(\pi x)$. We impose no-flux boundary conditions so that mass is conserved during evolution. In one-dimensional space, the $2$-Wasserstein distance between $\mu$ and $\nu$ is explicitly represented as 
\begin{equation}
    W_2^2(\mu,\nu) = \int_0^1 (U^{-1}(y)-V^{-1}(y))^2 dy
\end{equation}
where 
\begin{equation}
    U(x) = \int_{-\infty}^x d\mu, \quad V(x) = \int_{-\infty}^x d\nu.
\end{equation}
Using the formula we solve optimization problems by the $L^2$-gradient descent method with different time step sizes and demonstrate the accuracy by comparing the results with the exact solution (or, a numerical solution on a highly refined discretization as a proxy for the exact solution). We will refine the spatial grid as the time step size decreases so that the spatial discretization ensures sufficient accuracy. 

\begin{remark}
In the setting of Wasserstein gradient flows, the optimization problem that needs to be solved at each time step of scheme (\ref{eq:scheme}) can be challenging, especially in dimensions two and higher.
Indeed, not all of the coefficients $\gamma_{i,j}$ for either the second order scheme (\ref{eq:2ndScheme}) or the third order scheme (\ref{eq:3rdScheme}) are positive, which suggests a possible lack of convexity in the conventional sense for large time steps, a difficulty shared by e.g. the BDF2 method of \cite{matthes_plazotta}.
Although the authors of \cite{matthes_plazotta} note a different notion of convexity, namely that along generalized geodesics, holds in the case of BDF2, it is unclear if this property affords any practical help in finding the minimizer efficiently.
In particular, the approach in \cite{carillo_craig_wang_wei} of adapting Benamou \& Brenier's  dynamic reformulation of optimal transportation \cite{benamou_brenier} to the modified optimal transport problem that arises in backward Euler or Crank-Nicholson schemes, may no longer be applicable.
\end{remark}

\subsection{Heat equation}

Our result below demonstrates clear evidence of the expected order of convergence (while solutions remain smooth) of the new schemes constructed in light of the new stability analysis developed in the present paper, at least in the specific context of Wasserstein gradient flows. But even in this specific context, we leave the question of well-posedness of the variational problems involved to future work. Indeed, while the existence of minimizers is immediately verified independent of the time step size, the presence of multiple movement limiter terms with different signs may make questions of uniqueness harder than in \cite{matthes_plazotta}. 

We consider the heat equation $u_t(x,t) = u_{xx}(x,t)$ which is the Wasserstein gradient flow for the energy functional 
\begin{equation}
    \mathcal{E}(u) = \int u\log(u) dx.
\end{equation}
Given initial and boundary conditions as specified above, the exact solution to the heat equation is $u(x,t)=\frac{1}{2}+\frac{1}{4}\cos(\pi x)e^{-\pi^2t}$. We generate numerical solutions using two schemes (\ref{eq:2ndScheme}) and (\ref{eq:3rdScheme}) at time $t=\frac{1}{16}$. Table \ref{table:heat} shows relative $L^2$ errors and convergence orders for different numbers of time steps. 

\begin{table}[h!]
    \begin{minipage}{.5\linewidth}
    \centering
    \begin{tabular}{c  c  c}
    \hline 
    \# time steps & $L^2$ error & Order\\
    \hline \hline
    4 & 1.27E-04 & - \\ 
    \hline
    6 & 5.56E-05 & 2.04 \\ 
    \hline
    8 & 3.11E-05 & 2.02 \\ 
    \hline
    12 & 1.37E-05& 2.02 \\ 
    \hline
    16 & 7.64E-06 & 2.03 \\
    \hline
    24 &  3.39E-06 & 2.00 \\
    \hline
    \end{tabular}
    \end{minipage}
    \begin{minipage}{.5\linewidth}
    \centering
    \begin{tabular}{c  c  c}
    \hline 
    \# time steps & $L^2$ error & Order\\
    \hline \hline
    4 & 9.28E-06 & - \\ 
    \hline
    6 & 2.63E-06 & 3.11 \\ 
    \hline
    8 & 1.08E-06 & 3.09 \\ 
    \hline
    12 & 3.18E-07& 3.02 \\ 
    \hline
    16 & 1.35E-07 & 2.98 \\
    \hline
    24 &  4.19E-08 & 2.89 \\
    \hline
    \end{tabular}
    \end{minipage}
   \caption{Relative $L^2$ errors and corresponding orders of accuracy at time $t=\frac{1}{16}$ with increasing numbers of time steps when (left) the energy stable second order scheme (\ref{eq:2ndScheme}) and (right) the energy bounded third  order scheme (\ref{eq:3rdScheme}) are applied to the heat equation.} 
   \label{table:heat}
\end{table}

\begin{figure}[h!]
    \centering
       \includegraphics[clip, trim=2cm 6.5cm 2.5cm 7.5cm, width=0.7\textwidth]{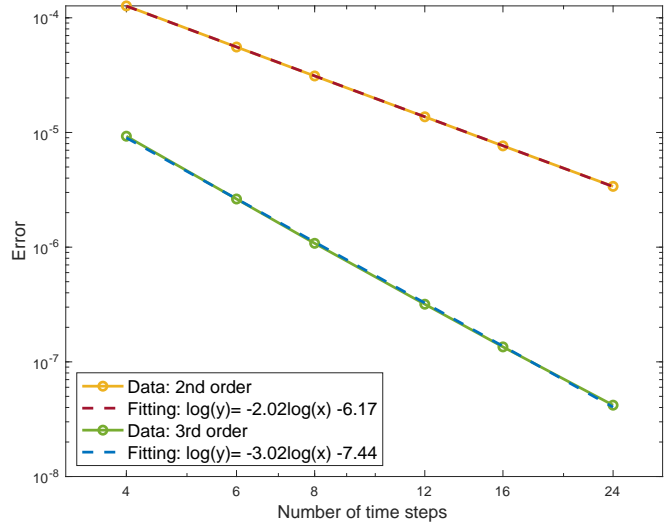}
    \caption{Errors in Table \ref{table:heat} on a log-log scale with line fitting. The magnitudes of slopes, $2.02$ and $3.02$,  represent the accuracy of each scheme.}
\end{figure}

\subsection{Porous medium equation}
We consider the porous medium equation $u_t(x,t) = (u^3(x,t))_{xx}$
which is the gradient flow for the energy
\begin{equation}
    \mathcal{E}(u) = \frac{1}{2}\int u^3 dx
\end{equation}
with respect to the quadratic Wasserstein metric. We apply second order and third  order schemes with various numbers of time steps and obtain numerical solutions at time $t=\frac{1}{8}$. We find a proxy for the exact solution by choosing an extremely small time step size so that the difference with another proxy generated by a smaller time step is less than $10^{-9}$. Table \ref{table:PME} tabulates the relative $L^2$ errors and orders of accuracy. 

\begin{table}[h!]
    \begin{minipage}{.5\linewidth}
    \centering
    \begin{tabular}{c  c  c}
    \hline 
    \# time steps & $L^2$ error & Order\\
    \hline \hline
    4 & 1.78E-04 & - \\ 
    \hline
    6 & 7.88E-05 & 2.01 \\ 
    \hline
    8 & 4.41E-05 & 2.02 \\ 
    \hline
    12 & 1.96E-05& 2.00 \\ 
    \hline
    16 & 1.10E-06 & 2.01 \\
    \hline
    24 &  4.91E-06 & 1.99 \\
    \hline
    32 & 2.77E-06 & 1.99 \\
    \hline
    \end{tabular}
    \end{minipage}
    \begin{minipage}{.5\linewidth}
    \centering
    \begin{tabular}{c  c  c}
    \hline 
    \# time steps & $L^2$ error & Order\\
    \hline \hline
    4 & 4.79E-05 & - \\ 
    \hline
    6 & 1.39E-05 & 3.05 \\ 
    \hline
    8 & 5.89E-06 & 2.98 \\ 
    \hline
    12 & 1.76E-06& 2.98 \\ 
    \hline
    16 & 7.52E-07 & 2.96 \\
    \hline
    24 &  2.29E-07 & 2.93 \\
    \hline
    32 & 9.95E-08 & 2.90\\
    \hline
    \end{tabular}
    \end{minipage}
   \caption{Relative $L^2$ errors and corresponding orders of accuracy  at time $t=\frac{1}{8}$ with increasing numbers of time steps when (left) the energy stable second order scheme (\ref{eq:2ndScheme}) and (right) the energy bounded third  order scheme (\ref{eq:3rdScheme}) are applied on the porous medium equation.}
   \label{table:PME}
\end{table}

\begin{figure}[h!]
    \centering
        \includegraphics[clip, trim=2cm 6.5cm 2.5cm 7.5cm, width=0.7\textwidth]{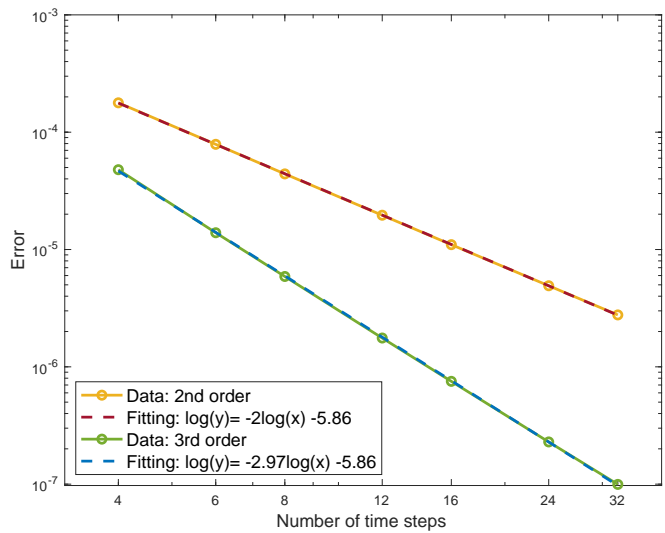}
    \caption{Errors in Table \ref{table:PME} on a log-log scale with line fitting. The magnitude of slope, $2.00$ and $2.97$,  represents the accuracy of each scheme.}
\end{figure}

\subsection{Nonlinear Fokker-Planck equation}
Next, we consider the nonlinear Fokker-Planck equation $ u_t(x,t) = (u(x,t)V_x(x))_x + (u^3(x,t))_{xx}$ where $V(x)=2+\cos(\pi x)$. This is the Wasserstein gradient flow of the energy 
\begin{equation}
    \mathcal{E}(u)=\int \left(u V + \frac{1}{2} u^3\right) dx.
\end{equation}
We find a highly accurate proxy for the exact solution at $t=\frac{1}{8}$ as in the previous example, and compare it with approximations which are obtained from our two schemes. Table \ref{table:FP} shows the $L^2$ errors and the orders of accuracy.

\begin{table}[h!]
    \begin{minipage}{.5\linewidth}
    \centering
    \begin{tabular}{c  c  c}
    \hline 
    \# time steps & $L^2$ error & Order\\
    \hline \hline
    6 & 9.09E-04 & -\\ 
    \hline
    8 & 5.04E-04 & 2.05 \\ 
    \hline
    12 & 2.21E-04& 2.03 \\ 
    \hline
    16 & 1.24E-04 & 2.01 \\
    \hline
    24 &  5.47E-05 & 2.02 \\
    \hline
    32 & 3.07E-05 & 2.01 \\
    \hline
    48 & 1.36E-05 & 2.01 \\
    \hline
    \end{tabular}
    \end{minipage}
    \begin{minipage}{.5\linewidth}
    \centering
    \begin{tabular}{c  c  c}
    \hline 
    \# time steps & $L^2$ error & Order\\
    \hline \hline
    8 & 4.30E-05 & - \\ 
    \hline
    12 & 1.24E-05 & 3.07 \\ 
    \hline
    16 & 5.21E-06 & 3.01 \\ 
    \hline
    24 & 1.58E-06& 2.94 \\ 
    \hline
    32 & 6.72E-07 & 2.97 \\
    \hline
    48 &  2.03E-07 & 2.95 \\
    \hline
    64 & 9.78E-08 & 2.54\\
    \hline
    \end{tabular}
    \end{minipage}
   \caption{Relative $L^2$ errors and corresponding orders of accuracy  at time $t=\frac{1}{8}$ with increasing numbers of time steps when (left) the energy stable second order scheme (\ref{eq:2ndScheme}) and (right) the energy bounded third  order scheme (\ref{eq:3rdScheme}) are applied to the nonlinear Fokker-Planck equation.}
   \label{table:FP}
\end{table}

\begin{figure}[h!]
    \centering
        \includegraphics[clip, trim=2cm 6.5cm 2.5cm 7.5cm, width=0.7\textwidth]{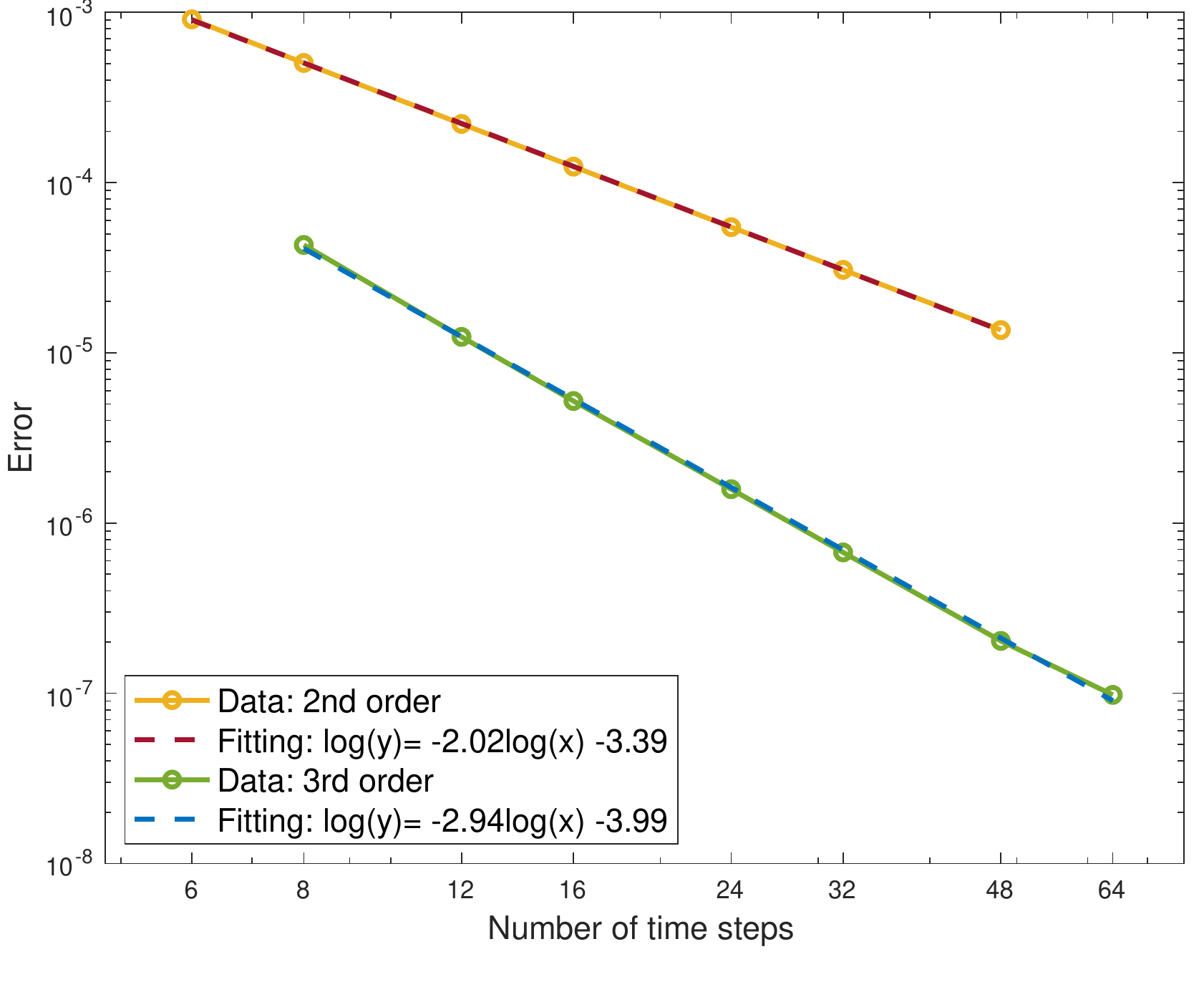}
    \caption{Errors in Table \ref{table:FP} on a log-log scale with line fitting. The magnitudes of slopes, $2.02$ and $2.94$,  represent the accuracy of each scheme.}
\end{figure}

\section{Conclusion}  In this paper we have established criteria to ensure energy dissipation and energy boundedness of multi-step, multi-stage, and mixed variational schemes in metric spaces from the viewpoint that certain finite graphs associated with the given scheme can be embedded into Euclidean spaces. We have worked through the process of applying the criteria on particular classes of schemes (multi-stage / two-step multi-stage diagonal schemes) and exhibited concrete examples of second and third order schemes that are unconditionally energy stable. The methods and their advertised accuracy are demonstrated on various examples in the 2-Wasserstein space.

\section{Appendix} 
\label{sec:Appendix}
\subsection{Proof of Claim \ref{claim:consistency}}
In this section, we assume $X$ to be a subspace of a finite dimensional Euclidean space $\mathbb{R}^n$ and derive second and third order consistency conditions for scheme (\ref{eq:scheme}) proposed in Claim \ref{claim:consistency}. From (\ref{eq:u_t}), we obtain higher order time derivatives of the exact solution to the gradient flow (\ref{eq:gradientFlow}) as follows:
\begin{equation}
\begin{dcases}
    u_t &= \Phi_1\\
    u_{tt}&=\Phi_2\\
    u_{ttt} &= \Phi_3+\Phi_4
\end{dcases}
\end{equation}
where 
\begin{equation}
    \begin{dcases}
        \Phi_1 &= -2(\nabla^2 D)^{-1}\nabla \mathcal{E}\\
        \Phi_2 &= -2(\nabla^2 D)^{-1}(\nabla^2 \mathcal{E}) \Phi_1\\
        \Phi_3 &= -2(\nabla^2 D)^{-1}(\nabla^2 \mathcal{E}) \Phi_2\\
        \Phi_4 &= -2(\nabla^2 D)^{-1}\begin{pmatrix}
         \Phi_1 ^T \nabla^2 \mathcal{E}_{x_1}\Phi_1\\
         \Phi_1 ^T \nabla^2 \mathcal{E}_{x_2}\Phi_1\\
         \vdots\\
         \Phi_1 ^T \nabla^2 \mathcal{E}_{x_n} \Phi_1\\
  \end{pmatrix}.
\end{dcases}
\end{equation}
We will show that the intermediate stages of the scheme (\ref{eq:scheme}) can be expanded in $k$ in terms of auxiliary functions $\Phi_1,\cdots,\Phi_4$ as:
\begin{equation}
    v_i = u_n+a_i\Phi_1 k +b_i\Phi_2 k^2+ (c_i\Phi_3+d_i\Phi_4)k^3 + O(k^4),
\end{equation}
for some $a_i,b_i,c_i,d_i\in\mathbb{R}$. 

Make an ansatz for the intermediate stage $v_i$ as 
\begin{equation}
    v_i = u_n + q_i^1k+q_i^2k^2+q_i^3k^3+O(k^4)
\end{equation}
for $q_i^1,q_i^2,q_i^3\in X$ and solve  Euler-Lagrange equation 
\begin{equation}
    \nabla \left\{\mathcal{E}(u)+\sum_j \frac{\gamma_{i,j}}{2k}D(u,v_j)\right\}=0.
\end{equation}
Apply straightforward Taylor expansion with respect to $k$ and consider each term (with respect to the order of $k$) separately. First, collect the $k^0$-terms. 
\begin{equation}
    \nabla \mathcal{E} + \sum_j \frac{\gamma_{i,j}}{2} (\nabla^2 D) (q_i^1-q_j^1) = 0.
\end{equation}
Since $q_j^1=a_j\Phi_1$ for all $j<0$ for some $a_j\in\mathbb{R}$, $q_i^1$ can be written as $q_i^1=a_i\Phi_1$ for $a_i\in\mathbb{R}$ as well, inductively. Then $a_i$ solves the equation
\begin{equation}
    1-\sum_j \gamma_{i,j}(a_i-a_j)=0,
\end{equation}
and this yields
\begin{equation}
    a_i=\frac{1}{\sum_j \gamma_{i,j}}\left(1+\sum_j \gamma_{i,j}a_j\right).
\end{equation}
Next, consider the $k^1$-terms. 
\begin{equation}
   \nabla^2 \mathcal{E} q_i^1 + \sum_j \frac{\gamma_{i,j}}{2} \nabla^2 D (q^2_i-q^2_j)= 0.
\end{equation}
Similarly, as $q_j^2 = b_j \Phi_2$ for $b_j\in\mathbb{R}$, $b_i$ solves the equation 
\begin{equation}
    a_i-\sum_{j} \gamma_{i,j}(b_i-b_j)=0,
\end{equation}
and hence
\begin{equation}
    b_i=\frac{1}{\sum_j \gamma_{i,j}}\left(a_i+\sum_j \gamma_{i,j}b_j\right).
\end{equation}
Finally, consider the $k^2$-terms. 
\begin{equation}
    (\nabla^2 \mathcal{E})q^2_i + \frac{1}{2} \begin{pmatrix}
     (q^1_i)^T (\nabla \mathcal{E}_{x_1}) (q^1_i)\\
     (q^1_i)^T (\nabla \mathcal{E}_{x_2}) (q^1_i)\\
     \vdots\\
     (q^1_i)^T (\nabla \mathcal{E}_{x_n}) (q^1_i)\\
    \end{pmatrix} 
    +\sum_{j} \frac{\gamma_{i,j}}{2} (\nabla^2 D)(q^3_i-q^3_j) =0.
\end{equation}
By setting $q_j^3 = c_j \Phi_3 + d_j \Phi_4$ for $c_j,d_j\in\mathbb{R}$, we get
\begin{equation}
    \begin{dcases}
        b_i-\sum_j \gamma_{i,j} (c_i-c_j)=0\\
        \frac{a_i^2}{2}-\sum_{j} \gamma_{i,j} (d_i-d_j)=0,
    \end{dcases}
\end{equation}
which is equivalent to
\begin{equation}
    \begin{dcases}
    c_i &= \frac{1}{\sum_j \gamma_{i,j}} \left(b_i+\sum_j \gamma_{i,j}c_j\right)\\
    d_i &= \sum_{j}\gamma_{i,j} \left(\frac{a_i^2}{2}+\sum_j\gamma_{i,j}d_j\right).
    \end{dcases}
\end{equation}

\subsection{Third  order accurate scheme}

In Section \ref{sec:consistency}, we presented a $2$-step $7$-stage diagonal scheme which is energy bounded and  third  order accurate. Here, we record the exact rational numbers of $\gamma_{i,j}$'s: 
\begin{equation}
    \begin{pmatrix}
     \gamma_{1,-1} & \gamma_{1,0} & 0 \\
    \gamma_{2,-1} & \gamma_{2,0} & \gamma_{2,1} \\
    \vdots\\
    \gamma_{6,-1} & \gamma_{6,0} & \gamma_{6,5} \\
    \gamma_{7,-1} & \gamma_{7,0} & \gamma_{7,6} \\
    \end{pmatrix} = \begin{pmatrix}
1/5& 324/25& 0\\
-67/100& 16/25& 249/20\\
-1/100& -19/25 & 1327/100 \\
  13/50& -71/50& 897/100\\
  1/20& -31/50& 69/10\\
 \gamma_{6,-1}&\gamma_{6,0}&\gamma_{6,5}\\
 \gamma_{7,-1}&\gamma_{7,0}&\gamma_{7,6}
\end{pmatrix}
\end{equation}
where
\begin{align*}
    &\gamma_{6,-1}=\frac{6738642394659375271309286924642199204}{
  499724271717869165338634999114429476375}\\
  &\gamma_{6,0}=-\frac{
   1490348725590513376673846530372322969031}{
   999448543435738330677269998228858952750}\\
  &\gamma_{6,5}=\frac{33204424381521663791982510017718750000}{
  3997794173742953322709079992915435811}\\
  &\gamma_{7,-1}=\frac{12657604782253956245795836543983271244969029}{
  68341222729403241230150248553811869282112250}\\
  &\gamma_{7,0}=-\frac{
   20148945983758481800702871507047428317759489}{
   34170611364701620615075124276905934641056125}\\
  &\gamma_{7,6}=\frac{384415962327102116281943490933129440840735787}{
  34170611364701620615075124276905934641056125}
  \end{align*}

\section{Acknowledgement}
Saem Han and Selim Esedo\=glu were supported by the National Science Foundation (NSF) grant DMS-2012015. Krishna Garikipati was supported by Toyota Research Institute grant \#849910. 

\bibliographystyle{plain}
\bibliography{references}

\end{document}